\documentclass[11pt]{article}
\usepackage[russian]{babel}
\usepackage[cp1251]{inputenc}
\usepackage{amsmath,latexsym}
\usepackage[psamsfonts]{amssymb}
\usepackage[dvips]{graphicx}
\usepackage[mathcal]{euscript}
\begin{document}

MSC2010.35G31

\begin{center}
\textbf{ Obtaining a representation of the solution to the Cauchy problem for \\ equation
high-order fractional derivative, \\
by the method of finding self-similar solutions}\\
\textbf{B.Yu.Irgashev}\\
\emph{Namangan  Engineering Construction Institute, \\ Namangan Branch of the Institute of Mathematics of the Republic of Uzbekistan.}\\
  \emph{E-mail: bahromirgasev@gmail.com}
\end{center}

\textbf{Abstract.} \emph{ In this paper, with the help of previously constructed self-similar solutions, a solution of the Cauchy problem for an equation of even order with a fractional Riemann-Liouville derivative of order $1<\alpha<2$ is obtained. }


\begin{center}
\textbf{1.Construction of self-similar solutions
}\end{center}
Consider the equation
\[L\left[ u \right] \equiv D_{0y}^\alpha u\left( {x,y} \right) - dD_{0x}^\beta u\left( {x,y} \right) = 0,\eqno(1)\]
here
\[q,p\in N,\,q - 1 < \alpha  \le q,\,\,p - 1 < \beta  \le p,\,\, q<p, \,d =  \pm 1,\]
and
$D_{0y}^\alpha ,\,D_{0x}^\beta  - $ Riemann-Liouville fractional differentiation operators, respectively, of orders $\alpha ,\beta $ :
\[D_{0y}^\alpha u\left( {x,y} \right) = \frac{1}{{\Gamma \left( {q - \alpha } \right)}}\frac{{{d^q}}}{{d{y^q}}}\int\limits_0^y {\frac{{u\left( {x,\tau } \right)d\tau }}{{{{\left( {y - \tau } \right)}^{\alpha  - q + 1}}}},\,y > 0,\,} \]
\[D_{0x}^\beta u\left( {x,y} \right) = \frac{1}{{\Gamma \left( {p - \beta } \right)}}\frac{{{d^p}}}{{d{x^p}}}\int\limits_0^x {\frac{{u\left( {\tau ,y} \right)d\tau }}{{{{\left( {x - \tau } \right)}^{\beta  - p + 1}}}},\,x > 0.\,} \]
Note that self-similar solutions of an equation with the usual derivative of the form
\[\frac{{{\partial ^p}u\left( {x,y} \right)}}{{\partial {x^p}}} - \frac{{{\partial ^q}u\left( {x,y} \right)}}{{\partial {y^q}}} = 0,\,p < q,\]
were constructed in [1].\\

In [2], using the methods of special operators, self-similar solutions of equation (1) were found in the cases $\frac{\beta }{2} \le \alpha  < \beta  \le 2,\,\,n - 1 < \beta  \le n \in N$.

In this paper, we construct self-similar solutions by a method that does not require knowledge of special operators.

We seek self-similar solutions to equation (1) in the form of the series
\[u\left( {x,y} \right) = {y^b}\sum\limits_{n = 0}^\infty  {{c_n}{{\left( {{x^a}{y^c}} \right)}^{n + \gamma }}}  = \sum\limits_{n = 0}^\infty  {{c_n}{x^{an + a\gamma }}{y^{cn + c\gamma  + b}}} ,\eqno(2)\]
here the parameters $a,b,c,\gamma  $ are to be defined.\\
We first introduce the notation
\[{\overline {\left( a \right)} _s} = a\left( {a - 1} \right)...\left( {a - \left( {s - 1} \right)} \right),\,\,{\overline {\left( a \right)} _0} = 1,\,\,{\overline {\left( a \right)} _1} = a.\]
Substitute (2) into (1). Then we will have
\[D_{0y}^\alpha u\left( {x,y} \right) = \frac{1}{{\Gamma \left( {q - \alpha } \right)}}\sum\limits_{n = 0}^\infty  {{c_n}{x^{an + a\gamma }}} \frac{{{d^q}}}{{d{y^q}}}\int\limits_0^y {\frac{{{\tau ^{cn + c\gamma  + b}}d\tau }}{{{{\left( {y - \tau } \right)}^{\alpha  - q + 1}}}}\,}  = \]
$$= \frac{1}{{\Gamma \left( {q - \alpha } \right)}}\sum\limits_{n = 0}^\infty  {{c_n}{x^{an + a\gamma }}} \frac{{{d^q}}}{{d{y^q}}}\int\limits_0^1 {\frac{{{y^{cn + c\gamma  + b + 1}}{z^{cn + c\gamma  + b}}dz}}{{{y^{\alpha  - q + 1}}{{\left( {1 - z} \right)}^{\alpha  - q + 1}}}}\,} = $$
\[ = \sum\limits_{n = 0}^\infty  {{{\overline {\left( {cn + c\gamma  + b - \alpha  + q} \right)} }_q}{c_n}{{\left( {{x^a}{y^c}} \right)}^{n + \gamma }}{y^{b - \alpha }}} \frac{{\Gamma \left( {cn + c\gamma  + b + 1} \right)}}{{\Gamma \left( {cn + c\gamma  + b - \alpha  + q + 1} \right)}},\eqno(3)\]
similarly
\[D_{0x}^\beta u\left( {x,y} \right) = \frac{1}{{\Gamma \left( {p - \beta } \right)}}\sum\limits_{n = 0}^\infty  {{c_n}{y^{cn + c\gamma  + b}}} \frac{{{d^p}}}{{d{x^p}}}\int\limits_0^x {\frac{{{\tau ^{an + a\gamma }}d\tau }}{{{{\left( {x - \tau } \right)}^{\beta  - p + 1}}}}\,}  = \]
\[ = \frac{1}{{\Gamma \left( {p - \beta } \right)}}\sum\limits_{n = 0}^\infty  {{c_n}{y^{cn + c\gamma  + b}}} \frac{{{d^p}}}{{d{x^p}}}\int\limits_0^x {\frac{{{x^{an + a\gamma  + 1}}{z^{an + a\gamma }}dz}}{{{x^{\beta  - p + 1}}{{\left( {1 - z} \right)}^{\beta  - p + 1}}}}\,}  = \]
\[ = \sum\limits_{n = 0}^\infty  {{{\overline {\left( {an + a\gamma  + p - \beta } \right)} }_p}{c_n}{{\left( {{x^a}{y^c}} \right)}^{n + \gamma }}{x^{ - \beta }}{y^b}} \frac{{\Gamma \left( {an + a\gamma  + 1} \right)}}{{\Gamma \left( {an + a\gamma  + 1 + p - \beta } \right)}}.\eqno(4)\]
Substituting (3), (4) into (1), we obtain
\[\sum\limits_{n = 0}^\infty  {{{\overline {\left( {cn + c\gamma  + b - \alpha  + q} \right)} }_q}{c_n}{{\left( {{x^a}{y^c}} \right)}^n}{x^\beta }{y^{ - \alpha }}} \frac{{\Gamma \left( {cn + c\gamma  + b + 1} \right)}}{{\Gamma \left( {cn + c\gamma  + b - \alpha  + q + 1} \right)}} = \]
\[ = d\sum\limits_{n = 0}^\infty  {{{\overline {\left( {an + a\gamma  + p - \beta } \right)} }_p}{c_n}{{\left( {{x^a}{y^c}} \right)}^n}} \frac{{\Gamma \left( {an + a\gamma  + 1} \right)}}{{\Gamma \left( {an + a\gamma  + 1 + p - \beta } \right)}} \Rightarrow \]
\[a = \beta ,\,c =  - \alpha  \Rightarrow \]
\[\sum\limits_{n = 0}^\infty  {{c_n}{{\left( {{x^\beta }{y^{ - \alpha }}} \right)}^{n + 1}}} \frac{{\Gamma \left( { - \alpha \left( {n + \gamma } \right) + b + 1} \right)}}{{\Gamma \left( { - \alpha \left( {n + \gamma  + 1} \right) + b + 1} \right)}} = \]
\[ = d\sum\limits_{n = 0}^\infty  {{{\overline {\left( {\beta n + \beta \gamma  + p - \beta } \right)} }_p}{c_n}{{\left( {{x^a}{y^c}} \right)}^n}} \frac{{\Gamma \left( {an + a\gamma  + 1} \right)}}{{\Gamma \left( {\beta n + \beta \gamma  + 1 + p - \beta } \right)}} \Rightarrow \]
got conditions for the parameter $\gamma :$
\[{\gamma _j} = 1 - \frac{j}{\beta },\,j = 1,2,...,p.\]
 Substituting the found parameters into representation (2), we find the formula for the coefficients ${c_n}$. We have
\[{c_n}\frac{{\Gamma \left( {\beta n + \beta \gamma  + 1} \right)}}{{\Gamma \left( {\beta n + \beta \gamma  + 1 - \beta } \right)}} = d{c_{n - 1}}\frac{{\Gamma \left( { - \alpha \left( {n + \gamma  - 1} \right) + b + 1} \right)}}{{\Gamma \left( { - \alpha \left( {n + \gamma } \right) + b + 1} \right)}} \Rightarrow \]
\[{c_n} = d{c_{n - 1}}\frac{{\Gamma \left( { - \alpha \left( {n + \gamma  - 1} \right) + b + 1} \right)\Gamma \left( {\beta n + \beta \gamma  + 1 - \beta } \right)}}{{\Gamma \left( { - \alpha \left( {n + \gamma } \right) + b + 1} \right)\Gamma \left( {\beta n + \beta \gamma  + 1} \right)}} = \]
\[ = {d^n}{c_0}\frac{{\Gamma \left( { - \alpha \gamma  + b + 1} \right)\Gamma \left( {\beta \gamma  + 1} \right)}}{{\Gamma \left( { - \alpha \left( {n + \gamma } \right) + b + 1} \right)\Gamma \left( {\beta \left( {n + \gamma } \right) + 1} \right)}},\]
if now
\[{c_0} = \frac{1}{{\Gamma \left( { - \alpha \gamma  + b + 1} \right)\Gamma \left( {\beta \gamma  + 1} \right)}},\]
then
\[{c_n} = \frac{{{d^n}}}{{\Gamma \left( { - \alpha \left( {n + 1 - \frac{j}{\beta }} \right) + b + 1} \right)\Gamma \left( {\beta \left( {n + 1 - \frac{j}{\beta }} \right) + 1} \right)}} = \]
\[ = \frac{{{d^n}}}{{\Gamma \left( { - \alpha n - \alpha  + \frac{\alpha }{\beta }j + b + 1} \right)\Gamma \left( {\beta n + \beta  - j + 1} \right)}}.\]
So we got the following self-similar solutions to equation (1):
\[{u_j}\left( {x,y} \right) = {y^b}{\left( {{x^\beta }{y^{ - \alpha }}} \right)^{1 - \frac{j}{\beta }}}\sum\limits_{n = 0}^\infty  {\frac{{{{\left( {d{x^\beta }{y^{ - \alpha }}} \right)}^n}}}{{\Gamma \left( { - \alpha n - \alpha  + \frac{\alpha }{\beta }j + b + 1} \right)\Gamma \left( {\beta n + \beta  - j + 1} \right)}}}  = \]
\[ = {y^b}{t^{1 - \frac{j}{\beta }}}{W_{\left( { - \alpha , - \alpha  + \frac{\alpha }{\beta }j + b + 1} \right),\left( {\beta ,\beta  - j + 1} \right)}}\left( {dt} \right),\,\,j = 1,2,...,p,\,\,\,t = {x^\beta }{y^{ - \alpha }},\eqno(5)\]
here
\[{W_{\left( {\mu ,a} \right),\left( {\nu ,b} \right)}}\left( z \right) = \sum\limits_{n = 0}^\infty  {\frac{{{z^n}}}{{\Gamma \left( {\mu n + a} \right)\Gamma \left( {\nu n + b} \right)}}} ,\,\mu ,\nu  \in R,\,a,b \in C,\,\mu  + \nu  > 0\]
- generalized Wright function [2].
Series (5) converges uniformly and it can be differentiated term by term, since $ \beta  - \alpha  > 0 $.
Note that there is no condition on the $ b $ parameter.Note that solutions of the form (5) coincide exactly with the solutions obtained in [2]. \\
Now let $\beta  = p \in N$, then the solutions of equation (1) are expressions of the form:
\[{u_s}\left( {x,y} \right) = {y^b}\sum\limits_{n = 0}^\infty  {\frac{{{d^n}{{\left( {x{y^{ - \frac{\alpha }{p}}}} \right)}^{pn + s}}}}{{\Gamma \left( { - \alpha n - \frac{\alpha }{p}s + b + 1} \right)\left( {pn + s} \right)!}}} ,\,\,s = 0,1,...,p - 1.\]
Then their linear combination is also a solution to equation (1)
\[u\left( {x,y} \right) = {y^b}\left( {{c_0}\sum\limits_{n = 0}^\infty  {\frac{{{d^n}{t^{pn}}}}{{\left( {pn} \right)!\Gamma \left( { - \alpha n + b + 1} \right)}}}  + {c_1}\sum\limits_{n = 0}^\infty  {\frac{{{d^n}{t^{pn + 1}}}}{{\left( {pn + 1} \right)!\Gamma \left( { - \alpha \left( {n + \frac{1}{p}} \right) + b + 1} \right)}} + ...} } \right.\]
\[ + \left. {{c_{p - 1}}\sum\limits_{n = 0}^\infty  {\frac{{{d^n}{t^{pn + p - 1}}}}{{\left( {p\left( {n + 1} \right) - 1} \right)!\Gamma \left( { - \alpha \left( {n + \frac{{p - 1}}{p}} \right) + b + 1} \right)}}} } \right),\eqno(6)\]
 where ${c_i} -$ arbitrary real numbers, $i=0,1,...,p-1 .$\\
Now let ${c_i} $, be such that ${c_i} \ne {c_j}$, for $i \ne j$ and $c_i^p = d$. Then from (6) we obtain
\[u\left( {x,y} \right) = {y^b}\phi \left( { - \frac{\alpha }{p},b + 1,cx{y^{ - \frac{\alpha }{p}}}} \right),\,{c^p} = 1,\eqno(7)\]
here
\[\phi \left( { - \delta ,\varepsilon ,z} \right) = \sum\limits_{k = 0}^\infty  {\frac{{{z^k}}}{{k!\Gamma \left( { - \delta k + \varepsilon } \right)}}} \]
- Wright function.\\
Consider the function
\[U\left( {x - \xi ,y - \eta } \right) = {\left( {y - \eta } \right)^b}\phi \left( { - \frac{\alpha }{p},b + 1,c\left( {x - \xi } \right){{\left( {y - \eta } \right)}^{ - \frac{\alpha }{p}}}} \right),\,y > \eta .\eqno(8)\]
Let us note some properties of function (8) that are obtained by direct computation.

\textbf{Лемма 1.}

1. $D_{0y}^\gamma \left( {{y^b}\phi \left( { - \frac{\alpha }{p},b + 1,cx{y^{ - \frac{\alpha }{p}}}} \right)} \right) = {y^{b - \gamma }}\phi \left( { - \frac{\alpha }{p},b + 1 - \gamma ,cx{y^{ - \frac{\alpha }{p}}}} \right),\,\gamma  \in R,$

2. $ D_{y\eta }^\gamma \left( {{{\left( {y - \eta } \right)}^b}\phi \left( { - \frac{\alpha }{p},b + 1,c\left( {x - \xi } \right){{\left( {y - \eta } \right)}^{ - \frac{\alpha }{p}}}} \right)} \right) = $

 $= {\left( {y - \eta } \right)^{b - \gamma }}\phi \left( { - \frac{\alpha }{p},b + 1 - \gamma ,c\left( {x - \xi } \right){{\left( {y - \eta } \right)}^{ - \frac{\alpha }{p}}}} \right),\,\gamma  \in R,$

3. $ \frac{{{\partial ^\gamma }}}{{\partial {x^\gamma }}}\left( {{y^b}\phi \left( { - \frac{\alpha }{p},b + 1,c\left( {x - \xi } \right){y^{ - \frac{\alpha }{p}}}} \right)} \right) = $

 $= {c^\gamma }{y^{b - \frac{\alpha }{p}\gamma }}\phi \left( { - \frac{\alpha }{p},b + 1 - \frac{\alpha }{p}\gamma ,c\left( {x - \xi } \right){y^{ - \frac{\alpha }{p}}}} \right),\,\gamma  \in \left\{ 0 \right\} \cup N,$

 4. $\left( {D_{\eta y}^\alpha  - d\frac{{{\partial ^p}}}{{\partial {x^p}}}} \right)D_{\eta y}^\gamma U\left( {x - \xi ,y - \eta } \right) = 0,$

 5. $ \left( {D_{y\eta }^\alpha  - {{\left( { - 1} \right)}^p}d\frac{{{\partial ^p}}}{{\partial {\xi ^p}}}} \right)D_{y\eta }^\gamma U\left( {x - \xi ,y - \eta } \right) = 0,$

 6. $ \left( {{{\left( { - 1} \right)}^q}D_{y\eta }^{\alpha  - q}\frac{{{\partial ^q}}}{{\partial {\eta ^q}}} - d{{\left( { - 1} \right)}^p}\frac{{{\partial ^p}}}{{\partial {\xi ^p}}}} \right)U\left( {x - \xi ,y - \eta } \right) = 0,\,q - 1 < \alpha  \le q.$

\begin{center}
\textbf{2.Cauchy problem
}\end{center}

In this section, using the found particular solutions (8), we obtain an explicit form of the solution to the Cauchy problem for an equation of the form:
\[D_{0y}^\alpha u\left( {x,y} \right) - {\left( { - 1} \right)^{n - 1}}\frac{{{\partial ^{2n}}u\left( {x,y} \right)}}{{\partial {x^{2n}}}} = f(x,y),\,\,1 < \alpha  < 2.\eqno(9)\]
Note that the Cauchy problem for higher order equations with a fractional derivative of order $0< \alpha <1$ was considered in [4], [5]. We will use ideas from these papers.

Following (7), we have
\[c_{1k}^{2n} = {\left( { - 1} \right)^{n - 1}} \Rightarrow {c_{1k}} = {e^{\frac{{n - 1 - 2k}}{{2n}}i\pi }},\,k = \overline {0,\left( {n - 1} \right)} ,\,{\mathop{\rm Re}\nolimits} \,{c_{1k}} > 0,\]
Consider the function
\[{\Gamma _b}\left( {x - \xi ,y - \eta } \right) = \left\{ \begin{array}{l}
\Gamma _b^1\left( {x - \xi ,y - \eta } \right),\,x > \xi ,\\
\Gamma _b^2\left( {x - \xi ,y - \eta } \right),\,x < \xi ,
\end{array} \right.\eqno(10)\]
where
\[\Gamma _b^1\left( {x - \xi ,y - \eta } \right) = \frac{{{{\left( {y - \eta } \right)}^b}}}{{2n}}\sum\limits_{k = 0}^{n - 1} {\left( { - {e^{\frac{{n - 1 - 2k}}{{2n}}i\pi }}} \right)\phi \left( { - \frac{\alpha }{{2n}},b + 1, - {e^{\frac{{n - 1 - 2k}}{{2n}}i\pi }}t} \right),\,} \]
\[\Gamma _b^2\left( {x - \xi ,y - \eta } \right) =  - \frac{{{{\left( {y - \eta } \right)}^b}}}{{2n}}\sum\limits_{k = 0}^{n - 1} {{e^{\frac{{n - 1 - 2k}}{{2n}}i\pi }}\phi \left( { - \frac{\alpha }{{2n}},b + 1, - {e^{\frac{{n - 1 - 2k}}{{2n}}i\pi }}\left( { - t} \right)} \right)}, \]
\[t = \frac{{\left( {x - \xi } \right)}}{{{{\left( {y - \eta } \right)}^{\frac{\alpha }{{2n}}}}}} .\]
We have
\[\frac{{{\partial ^s}\Gamma _b^1\left( {x - \xi ,y - \eta } \right)}}{{\partial {x^s}}} = \frac{1}{{2n}}\sum\limits_{k = 0}^{n - 1} {\left( { - {e^{\frac{{n - 1 - 2k}}{{2n}}i\pi }}} \right)\frac{{{\partial ^s}\left( {{{\left( {y - \eta } \right)}^b}\phi \left( { - \frac{\alpha }{{2n}},b + 1, - {e^{\frac{{n - 1 - 2k}}{{2n}}i\pi }}t} \right)} \right)}}{{\partial {x^s}}}}  = \]
\[ = \frac{1}{{2n}}\sum\limits_{k = 0}^{n - 1} {{{\left( { - {e^{\frac{{n - 1 - 2k}}{{2n}}i\pi }}} \right)}^{s + 1}}{{\left( {y - \eta } \right)}^{b - \frac{\alpha }{{2n}}s}}\phi \left( { - \frac{\alpha }{{2n}},b + 1 - \frac{\alpha }{{2n}}s, - {e^{\frac{{n - 1 - 2k}}{{2n}}i\pi }}t} \right)}  \Rightarrow \]
\[{\left. {\frac{{{\partial ^s}\Gamma _b^1\left( {x - \xi ,y - \eta } \right)}}{{\partial {x^s}}}} \right|_{x = \xi }} = \frac{{{{\left( {y - \eta } \right)}^{b - \frac{\alpha }{{2n}}s}}}}{{2n\Gamma \left( {b + 1 - \frac{\alpha s}{{2n}}} \right)}}\sum\limits_{k = 0}^{n - 1} {{{\left( { - {e^{\frac{{n - 1 - 2k}}{{2n}}i\pi }}} \right)}^{s + 1}}} .\]
Similarly
\[{\left. {\frac{{{\partial ^s}\Gamma _b^2\left( {x - \xi ,y - \eta } \right)}}{{\partial {x^s}}}} \right|_{x = \xi }} =  - \frac{{{{\left( {y - \eta } \right)}^{b - \frac{\alpha }{{2n}}s}}}}{{2n\Gamma \left( {b + 1 - \frac{\alpha s }{{2n}}} \right)}}\sum\limits_{k = 0}^{n - 1} {{{\left( {{e^{\frac{{n - 1 - 2k}}{{2n}}i\pi }}} \right)}^{s + 1}}} .\]
Hence we have if
\[s = \left( {2n - 1} \right)\,\left( {\bmod \,2n} \right),\]
then
\[{\left. {\frac{{{\partial ^s}\Gamma _b^1\left( {x - \xi ,y - \eta } \right)}}{{\partial {x^s}}}} \right|_{x = \xi }} - {\left. {\frac{{{\partial ^s}\Gamma _b^2\left( {x - \xi ,y - \eta } \right)}}{{\partial {x^s}}}} \right|_{x = \xi }} = {\left( { - 1} \right)^{n - 1}}\frac{{{{\left( {y - \eta } \right)}^{b - \frac{\alpha }{{2n}}s}}}}{{\Gamma \left( {b + 1 - \frac{\alpha s }{{2n}}} \right)}},\]
if
\[s \ne \left( {2n - 1} \right)\,\left( {\bmod \,2n} \right),\]
then
\[{\left. {\frac{{{\partial ^s}\Gamma _b^1\left( {x - \xi ,y - \eta } \right)}}{{\partial {x^s}}}} \right|_{x = \xi }} - {\left. {\frac{{{\partial ^s}\Gamma _b^2\left( {x - \xi ,y - \eta } \right)}}{{\partial {x^s}}}} \right|_{x = \xi }} = \]
\[ = \frac{{{{\left( {y - \eta } \right)}^{b - \frac{\alpha }{{2n}}s}}}}{{2n\Gamma \left( {b + 1 - \frac{\alpha }{{2n}}} \right)}}\left( {\sum\limits_{k = 0}^{n - 1} {{{\left( { - {e^{\frac{{n - 1 - 2k}}{{2n}}i\pi }}} \right)}^{s + 1}}}  + \sum\limits_{k = 0}^{n - 1} {{{\left( {{e^{\frac{{n - 1 - 2k}}{{2n}}i\pi }}} \right)}^{s + 1}}} } \right) = \]
\[ = \frac{{{{\left( {y - \eta } \right)}^{b - \frac{\alpha }{{2n}}s}}}}{{2n\Gamma \left( {b + 1 - \frac{\alpha }{{2n}}} \right)}}{e^{\frac{{n - 1}}{{2n}}\left( {s + 1} \right)i\pi }}\sum\limits_{k = 0}^{n - 1} {{e^{ - \frac{{k\left( {s + 1} \right)}}{n}i\pi }}\left( {{{\left( { - 1} \right)}^{s + 1}} + 1} \right)}  = 0.\]
So we have proved the following lemma.

\textbf{Lemma 2.} For $ s\in N $, the following relation holds:
\[{\left. {\frac{{{\partial ^s}\Gamma _b^1\left( {x - \xi ,y - \eta } \right)}}{{\partial {x^s}}}} \right|_{x = \xi }} - {\left. {\frac{{{\partial ^s}\Gamma _b^2\left( {x - \xi ,y - \eta } \right)}}{{\partial {x^s}}}} \right|_{x = \xi }} = \]
\[ = {\left( { - 1} \right)^{n - 1}}\frac{{{{\left( {y - \eta } \right)}^{b - \frac{\alpha }{{2n}}s}}}}{{\Gamma \left( {b + 1 - \frac{\alpha }{{2n}}} \right)}}\left\{ \begin{array}{l}
1,\,\,s = \left( {2n - 1} \right)\,\left( {\bmod \,2n} \right),\\
0,\,\,s \ne \left( {2n - 1} \right)\,\left( {\bmod \,2n} \right).
\end{array} \right.\]
The results of Lemma 2 coincide with the results of [5].
In what follows, we need the asymptotics of the Wright function for large values of the variable. The main results on asymptotics were obtained by Wright [see 3]. In particular, the following theorem is proved.

\textbf{Theorem 1} [see 3]. If $\left| {\arg \,y} \right| \le \min \left\{ {\frac{3}{2}\pi \left( {1 - \sigma } \right),\pi } \right\} - \varepsilon ,\,\,0 < \sigma  < 1,$ then
\[\phi \left( { - \sigma ,\beta ,t} \right) = {Y^{\frac{1}{2} - \beta }}{e^{ - Y}}\left\{ {\sum\limits_{m = 0}^{M - 1} {{A_m}{Y^{ - m}}}  + O\left( {{Y^{ - M}}} \right)} \right\},\,\,\left| t \right| \to \infty ,\eqno(11)\]
 here $ Y = \left( {1 - \sigma } \right){\left( {{\sigma ^\sigma }y} \right)^{\frac{1}{{1 - \sigma }}}},\,\,y =  - t,\,$ the coefficients of $ {A_m} $ depend on $\sigma ,\beta $.\\
If $n=3,4,...$, то $ 0 < \frac{\alpha }{{2n}} < \frac{1}{n} \le \frac{1}{3},$  and the relation
\[\left| {\arg \,y} \right| = \left| {\arg \,\left( { - t} \right)} \right| = \left| {\frac{{n - 1 - 2k}}{{2n}}\pi } \right| \le \pi  - \varepsilon , k=0,1,...,n-1.\eqno(12)\]
If $n=2$, then for $ 0 < \frac{\alpha }{4} \le \frac{1}{3},$ we have
\[\left| {\arg \,y} \right| = \left| {\arg \,\left( { - t} \right)} \right| = \left| { \pm \frac{\pi }{4}} \right| < \pi ,\eqno(13)\]
and at $ \frac{1}{3} < \frac{\alpha }{4} < \frac{1}{2},$ get
\[\left| {\arg \,y} \right| = \left| {\arg \,\left( { - t} \right)} \right| = \left| { \pm \frac{\pi }{4}} \right| < \frac{3}{2}\pi \left( {1 - \frac{\alpha }{4}} \right).\eqno(14)\]
Taking into account (12) - (14), we conclude that for $ n = 2,3, ... $ we always have relation (11). Now we write (10) in the form
\[{\Gamma _b}\left( {x - \xi ,y - \eta } \right) = \frac{{{{\left( {y - \eta } \right)}^b}}}{{2n}}\sum\limits_{k = 0}^{n - 1} {\left( { - {e^{\frac{{n - 1 - 2k}}{{2n}}i\pi }}} \right)\phi \left( { - \frac{\alpha }{{2n}},b + 1, - {e^{\frac{{n - 1 - 2k}}{{2n}}i\pi }}\left| t \right|} \right),\,}\eqno(15) \]
where
\[\left| t \right| = \frac{{\left| {x - \xi } \right|}}{{{{\left( {y - \eta } \right)}^{\frac{\alpha }{{2n}}}}}},\]
then, taking into account (11), (12) - (14), for large values of $\left| t \right|$, the estimate obtained in [5], [6] is valid for (15):
\[\left| {\phi \left( { - \frac{\alpha }{{2n}},b + 1, - {e^{\frac{{n - 1 - 2k}}{{2n}}i\pi }}\left| t \right|} \right)} \right| \le C{\left| t \right|^{\frac{{ - 2n}}{{2n - \alpha }}\left( {b + \frac{1}{2}} \right)}}\exp \left( { - \sigma {{\left| t \right|}^{\frac{{2n}}{{2n - \alpha }}}}} \right),\eqno(16)\]
here
\[\sigma  = \left( {1 - \frac{\alpha }{{2n}}} \right){\left( {\frac{\alpha }{{2n}}} \right)^{\frac{\alpha }{{2n - \alpha }}}}\cos \frac{{n - 1}}{{2n - \alpha }}\pi ,\,0 < C - const,\]
or applying Lemma 1 (properties 1,3), for $x \ne 0,\,\gamma  \in R,\,s \in N \cup \left\{ 0 \right\},$ we obtain
\[\left| {\frac{{{\partial ^s}}}{{\partial {x^s}}}\left( {D_{0y}^\gamma {\Gamma _b}\left( {x,y} \right)} \right)} \right| \le M{y^{b - \gamma  - \frac{\alpha }{{2n}}\left( {s - \theta } \right)}}{\left| x \right|^{ - \theta }}\exp \left( { - \sigma {{\left| x \right|}^{\frac{{2n}}{{2n - \alpha }}}}{y^{ - \frac{\alpha }{{2n - \alpha }}}}} \right),\eqno(17)\]
here
\[\left| x \right|{y^{ - \frac{\alpha }{{2n}}}} \to  + \infty, \]  \[\theta  = \frac{{2n}}{{2n - \alpha }}\left( {\frac{1}{2} + b - \gamma  - \frac{\alpha }{{2n}}s} \right).\]
Note that estimate (17) coincides with the estimate obtained in [5], for $b = \alpha  - 1 - \frac{\alpha }{{2n}}.$

\textbf{Lemma 3.} For $j,s\, \in N$ and any function  $h\left( x \right) \in C\left( R \right)$ such that
\[\left| {h\left( x \right)} \right| \le M\exp \left( {c{{\left| x \right|}^{\frac{{2n}}{{2n - \alpha }}}}} \right),\,\,c < \sigma ,\,\,\left| x \right| \to \infty,0<M-const,\eqno(18)\]
equality holds
\[\mathop {\lim }\limits_{y \to 0} \int\limits_{ - \infty }^{ + \infty } {h\left( \xi  \right)D_{0y}^{\alpha  - s}{\Gamma _{\alpha  - \frac{\alpha }{{2n}} - j}}\left( {x - \xi ,y} \right)d\xi }  = -h\left( x \right) \cdot \left\{ \begin{array}{l}
0,\,\,j \ne s,\\
1,\,\,\,j = s.
\end{array} \right.\]

\textbf{Proof.} First, we calculate the integral
\[\int\limits_{ - \infty }^{ + \infty } {D_{0y}^{\alpha  - s}{\Gamma _{\alpha  - \frac{\alpha }{{2n}} - j}}\left( {x - \xi ,y} \right)d\xi }  = \int\limits_{ - \infty }^x {D_{0y}^{\alpha  - s}{\Gamma _{\alpha  - \frac{\alpha }{{2n}} - j}}\left( {x - \xi ,y} \right)d\xi }  + \]
\[ + \int\limits_x^{ + \infty } {D_{0y}^{\alpha  - s}{\Gamma _{\alpha  - \frac{\alpha }{{2n}} - j}}\left( {x - \xi ,y} \right)d\xi }  = \]
\[\begin{array}{l}
 =  - \sum\limits_{k = 0}^{n - 1} {\frac{{{e^{\frac{{n - 1 - 2k}}{{2n}}\pi i}}}}{{2n}}} \int\limits_{ - \infty }^x {D_{0y}^{\alpha  - s}\left( {{y^{\alpha  - \frac{\alpha }{{2n}} - j}}\phi \left( { - \frac{\alpha }{{2n}},\alpha  - \frac{\alpha }{{2n}} - j + 1, - {e^{\frac{{n - 1 - 2k}}{{2n}}\pi i}}\left( {x - \xi } \right){y^{ - \frac{\alpha }{{2n}}}}} \right)} \right)d\xi }  - \\
 - \sum\limits_{k = 0}^{n - 1} {\frac{{{e^{\frac{{n - 1 - 2k}}{{2n}}\pi i}}}}{{2n}}} \int\limits_x^{ + \infty } {D_{0y}^{\alpha  - s}\left( {{y^{\alpha  - \frac{\alpha }{{2n}} - j}}\phi \left( { - \frac{\alpha }{{2n}},\alpha  - \frac{\alpha }{{2n}} - j + 1,{e^{\frac{{n - 1 - 2k}}{{2n}}\pi i}}\left( {x - \xi } \right){y^{ - \frac{\alpha }{{2n}}}}} \right)} \right)d\xi }  =
\end{array}\]
\[ =  - \sum\limits_{k = 0}^{n - 1} {\frac{{{e^{\frac{{n - 1 - 2k}}{{2n}}\pi i}}}}{{2n}}} \int\limits_{ - \infty }^x {{y^{s - \frac{\alpha }{{2n}} - j}}\phi \left( { - \frac{\alpha }{{2n}},s - \frac{\alpha }{{2n}} - j + 1, - {e^{\frac{{n - 1 - 2k}}{{2n}}\pi i}}\left( {x - \xi } \right){y^{ - \frac{\alpha }{{2n}}}}} \right)d\xi }  - \]
\[ - \sum\limits_{k = 0}^{n - 1} {\frac{{{e^{\frac{{n - 1 - 2k}}{{2n}}\pi i}}}}{{2n}}} \int\limits_x^{ + \infty } {{y^{s - \frac{\alpha }{{2n}} - j}}\phi \left( { - \frac{\alpha }{{2n}},s - \frac{\alpha }{{2n}} - j + 1,{e^{\frac{{n - 1 - 2k}}{{2n}}\pi i}}\left( {x - \xi } \right){y^{ - \frac{\alpha }{{2n}}}}} \right)d\xi }.\]
Now let's calculate each integral separately
\[ - {e^{\frac{{n - 1 - 2k}}{{2n}}\pi i}}\int\limits_{ - \infty }^x {{y^{s - \frac{\alpha }{{2n}} - j}}\phi \left( { - \frac{\alpha }{{2n}},s - \frac{\alpha }{{2n}} - j + 1, - {e^{\frac{{n - 1 - 2k}}{{2n}}\pi i}}\left( {x - \xi } \right){y^{ - \frac{\alpha }{{2n}}}}} \right)d\xi } \]
\[ =  - {e^{\frac{{n - 1 - 2k}}{{2n}}\pi i}}\sum\limits_{m = 0}^\infty  {{y^{s - \frac{\alpha }{{2n}} - j}}} \int\limits_{ - \infty }^x {\frac{{{{\left( { - {e^{\frac{{n - 1 - 2k}}{{2n}}\pi i}}{y^{ - \frac{\alpha }{{2n}}}}\left( {x - \xi } \right)} \right)}^m}}}{{m!\Gamma \left( { - \frac{\alpha }{{2n}}m + s - \frac{\alpha }{{2n}} - j + 1} \right)}}d\xi }  = \]
\[ =  - {y^{s - j}}\sum\limits_{m = 0}^\infty  {\left. {\frac{{{{\left( { - {e^{\frac{{n - 1 - 2k}}{{2n}}\pi i}}{y^{ - \frac{\alpha }{{2n}}}}\left( {x - \xi } \right)} \right)}^{m + 1}}}}{{\left( {m + 1} \right)!\Gamma \left( { - \frac{\alpha }{{2n}}m + s - \frac{\alpha }{{2n}} - j + 1} \right)}}} \right|_{\xi  = -\infty }^{\xi  = x}}  = \]
\[ =  - {y^{s - j}}\left. {\left( {\sum\limits_{m = 0}^\infty  {\frac{{{{\left( { - {e^{\frac{{n - 1 - 2k}}{{2n}}\pi i}}{y^{ - \frac{\alpha }{{2n}}}}\left( {x - \xi } \right)} \right)}^m}}}{{m!\Gamma \left( { - \frac{\alpha }{{2n}}m + s - j + 1} \right)}} - \frac{1}{{\Gamma \left( {s - j + 1} \right)}}} } \right)} \right|_{\xi  =- \infty }^{\xi  = x} = \]
\[ =  - {y^{s - j}}\left. {\left( {\phi \left( { - \frac{\alpha }{{2n}},s - j + 1, - {e^{\frac{{n - 1 - 2k}}{{2n}}\pi i}}{y^{ - \frac{\alpha }{{2n}}}}\left( {x - \xi } \right)} \right) - \frac{1}{{\Gamma \left( {s - j + 1} \right)}}} \right)} \right|_{\xi  =- \infty }^{\xi  = x} = \]
\[ = -\frac{{{y^{s - j}}}}{{\Gamma \left( {s - j + 1} \right)}}.\]
Similarly
\[ - {e^{\frac{{n - 1 - 2k}}{{2n}}\pi i}}\int\limits_x^{ + \infty } {{y^{s - \frac{\alpha }{{2n}} - j}}\phi \left( { - \frac{\alpha }{{2n}},s - \frac{\alpha }{{2n}} - j + 1,{e^{\frac{{n - 1 - 2k}}{{2n}}\pi i}}\left( {x - \xi } \right){y^{ - \frac{\alpha }{{2n}}}}} \right)d\xi }  = -\frac{{{y^{s - j}}}}{{\Gamma \left( {s - j + 1} \right)}}.\]
from here we get
\[\int\limits_{ - \infty }^{ + \infty } {D_{0y}^{\alpha  - s}{\Gamma _{\alpha  - \frac{\alpha }{{2n}} - j}}\left( {x - \xi ,y} \right)d\xi }  =- \frac{{{y^{s - j}}}}{{\Gamma \left( {s - j + 1} \right)}}.\eqno(19)\]
Note that for $ s <j $ we have
\[\int\limits_{ - \infty }^{ + \infty } {D_{0y}^{\alpha  - s}{\Gamma _{\alpha  - \frac{\alpha }{{2n}} - j}}\left( {x - \xi ,y} \right)d\xi }  = 0,\,\,s < j.\eqno(20)\]
Now
\[\int\limits_{ - \infty }^{ + \infty } {h\left( \xi  \right)D_{0y}^{\alpha  - s}{\Gamma _{\alpha  - \frac{\alpha }{{2n}} - j}}\left( {x - \xi ,y} \right)d\xi }  = \int\limits_{ - \infty }^x {h\left( \xi  \right)D_{0y}^{\alpha  - s}{\Gamma _{\alpha  - \frac{\alpha }{{2n}} - j}}\left( {x - \xi ,y} \right)d\xi }  + \]
\[ + \int\limits_x^{ + \infty } {h\left( \xi  \right)D_{0y}^{\alpha  - s}{\Gamma _{\alpha  - \frac{\alpha }{{2n}} - j}}\left( {x - \xi ,y} \right)d\xi }  = {I_1} + {I_2},\]
We calculate each term separately
\[{I_1} = \int\limits_{ - \infty }^x {h\left( \xi  \right)D_{0y}^{\alpha  - s}{\Gamma _{\alpha  - \frac{\alpha }{{2n}} - j}}\left( {x - \xi ,y} \right)d\xi }  = \]
\[ =  - \sum\limits_{k = 0}^{n - 1} {{e^{\frac{{n - 1 - 2k}}{{2n}}\pi i}}{y^{s - \frac{\alpha }{{2n}} - j}}} \int\limits_{ - \infty }^x {h\left( \xi  \right)\phi \left( { - \frac{\alpha }{{2n}},s - \frac{\alpha }{{2n}} - j + 1, - {e^{\frac{{n - 1 - 2k}}{{2n}}\pi i}}{y^{ - \frac{\alpha }{{2n}}}}\left( {x - \xi } \right)} \right)d\xi }  = \]
\[ = \left( {t = \left( {x - \xi } \right){y^{ - \frac{\alpha }{{2n}}}}} \right) = \]
\[ =  - \sum\limits_{k = 0}^{n - 1} {{e^{\frac{{n - 1 - 2k}}{{2n}}\pi i}}{y^{s - j}}} \int\limits_0^{ + \infty } {h\left( {x - t{y^{\frac{\alpha }{{2n}}}}} \right)\phi \left( { - \frac{\alpha }{{2n}},s - \frac{\alpha }{{2n}} - j + 1, - {e^{\frac{{n - 1 - 2k}}{{2n}}\pi i}}t} \right)d\xi }  = \]
\[ =  - \sum\limits_{k = 0}^{n - 1} {{e^{\frac{{n - 1 - 2k}}{{2n}}\pi i}}{y^{s - j}}} \int\limits_0^{ + \infty } {\left( {h\left( {x - t{y^{\frac{\alpha }{{2n}}}}} \right) - h\left( x \right)} \right)\phi \left( { - \frac{\alpha }{{2n}},s - \frac{\alpha }{{2n}} - j + 1, - {e^{\frac{{n - 1 - 2k}}{{2n}}\pi i}}t} \right)d\xi  + } \]
\[ + h\left( x \right)\left( { - \sum\limits_{k = 0}^{n - 1} {{e^{\frac{{n - 1 - 2k}}{{2n}}\pi i}}{y^{s - j}}} \int\limits_0^{ + \infty } {\phi \left( { - \frac{\alpha }{{2n}},s - \frac{\alpha }{{2n}} - j + 1, - {e^{\frac{{n - 1 - 2k}}{{2n}}\pi i}}t} \right)d\xi } } \right) = \]
\[ = {J_1} + {J_2},\]
Hence we have
\[\mathop {\lim }\limits_{y \to  + 0} {I_1} = \mathop {\lim }\limits_{y \to  + 0} {J_1} + \mathop {\lim }\limits_{y \to  + 0} {J_2}.\]
Let's find each limit separately. Using (19), (20), we have
\[\mathop {\lim }\limits_{y \to  + 0} {J_2} =- h\left( x \right)\mathop {\lim }\limits_{y \to  + 0} \frac{{{y^{s - j}}}}{{2\Gamma \left( {s - j + 1} \right)}} = \left\{ \begin{array}{l}
-\frac{1}{2}h\left( x \right),\,s = j,\\
0,\,\,s \ne j.
\end{array} \right.\]
Further, if $s \ne j,$ then $\mathop {\lim }\limits_{y \to  + 0} {J_1} = 0.$ Let $ s = j $, then
\[\mathop {\lim }\limits_{y \to  + 0} \int\limits_0^{ + \infty } {\left( {h\left( {x - t{y^{\frac{\alpha }{{2n}}}}} \right) - h\left( x \right)} \right)\phi \left( { - \frac{\alpha }{{2n}},s - \frac{\alpha }{{2n}} - j + 1, - {e^{\frac{{n - 1 - 2k}}{{2n}}\pi i}}t} \right)d\xi }  = \]
\[ = \mathop {\lim }\limits_{y \to  + 0} \int\limits_0^N {\left( {h\left( {x - t{y^{\frac{\alpha }{{2n}}}}} \right) - h\left( x \right)} \right)\phi \left( { - \frac{\alpha }{{2n}},s - \frac{\alpha }{{2n}} - j + 1, - {e^{\frac{{n - 1 - 2k}}{{2n}}\pi i}}t} \right)d\xi }  + \]
\[ + \mathop {\lim }\limits_{y \to  + 0} \int\limits_N^{ + \infty } {\left( {h\left( {x - t{y^{\frac{\alpha }{{2n}}}}} \right) - h\left( x \right)} \right)\phi \left( { - \frac{\alpha }{{2n}},s - \frac{\alpha }{{2n}} - j + 1, - {e^{\frac{{n - 1 - 2k}}{{2n}}\pi i}}t} \right)d\xi }  = 0,\]
here $ 0 <N $. The last equality holds due to the fact that for an arbitrary small number $\epsilon > 0 $, you can choose $ N $ so that the sublimit expressions will be less than $1/2\epsilon$ (the first expression, due to the continuity the function $ h (x) $, and the second is due to estimates (16) and (18)). \\
Means
\[\mathop {\lim }\limits_{y \to  + 0} {I_1} = \left\{ \begin{array}{l}
-\frac{1}{2}h\left( x \right),\,s = j,\\
0,\,\,s \ne j.
\end{array} \right.\]
It is shown similarly that
\[\mathop {\lim }\limits_{y \to  + 0} {I_2} = \left\{ \begin{array}{l}
-\frac{1}{2}h\left( x \right),\,s = j,\\
0,\,\,s \ne j.
\end{array} \right.\]
So finally we have
\[\mathop {\lim }\limits_{y \to 0} \int\limits_{ - \infty }^{ + \infty } {h\left( \xi  \right)D_{0y}^{\alpha  - s}{\Gamma _{\alpha  - \frac{\alpha }{{2n}} - j}}\left( {x - \xi ,y} \right)d\xi }  =- h\left( x \right) \cdot \left\{ \begin{array}{l}
0,\,\,j \ne s,\\
1,\,\,\,j = s.
\end{array} \right.\]
\textbf{Lemma 3 is proved.}

\textbf{Cauchy problem.}Find the solution $u\left( {x,y} \right)$ to equation (9)
in the region $ D = \left\{ {\left( {x,y} \right): - \infty  < x <  + \infty ,0 < y} \right\},$ satisfying the following conditions:

1) $D_{0y}^\alpha u\left( {x,y} \right),\frac{{{\partial ^{2n}}u\left( {x,y} \right)}}{{\partial {x^{2n}}}} \in C\left( D \right),D_{0y}^{\alpha  - 1}u\left( {x,y} \right),D_{0y}^{\alpha  - 2}u\left( {x,y} \right) \in C\left( {\overline D } \right);$

2)$\mathop {\lim }\limits_{y \to  + 0} D_{oy}^{\alpha  - 1}u\left( {x,y} \right) = \varphi \left( x \right),\,\mathop {\lim }\limits_{y \to  + 0} D_{oy}^{\alpha  - 2}u\left( {x,y} \right) = \psi \left( x \right).$\\
The specified functions satisfy the constraints:
\[f\left( {x,y} \right) \in C\left( {\overline D } \right),\varphi \left( x \right),\psi \left( x \right) \in C\left( R \right),\]
\[\left| {\varphi \left( x \right)} \right| \le M\exp \left( {k{{\left| x \right|}^{\frac{{2n}}{{2n - \alpha }}}}} \right),\,\left| x \right| \to  + \infty ,\]
\[\left| {\psi \left( x \right)} \right| \le M\exp \left( {k{{\left| x \right|}^{\frac{{2n}}{{2n - \alpha }}}}} \right),\,\left| x \right| \to  + \infty ,\]
\[\left| {f\left( {x,y} \right)} \right| \le M\exp \left( {k{{\left| x \right|}^{\frac{{2n}}{{2n - \alpha }}}}} \right),\,\left| x \right| \to  + \infty ,$$
$$k < \sigma, 0<M-constant.\]

\textbf{Theorem 2.} The solution to the Cauchy problem has the form
\[u\left( {x,y} \right) =- \int\limits_{ - \infty }^{ + \infty } {\varphi \left( \xi  \right){\Gamma _{\alpha  - \frac{\alpha }{{2n}} - 1}}\left( {x - \xi ,y} \right)d\xi }  - \int\limits_{ - \infty }^{ + \infty } {\psi \left( \xi  \right){\Gamma _{\alpha  - \frac{\alpha }{{2n}} - 2}}\left( {x - \xi ,y} \right)d\xi }  - \]
\[ - \int\limits_0^y {\int\limits_{ - \infty }^{ + \infty } {f\left( {\xi ,\eta } \right){\Gamma _{\alpha  - \frac{\alpha }{{2n}} - 1}}\left( {x - \xi ,y - \eta } \right)} } d\xi d\eta .\]

\textbf{Proof.} It is easy to verify that the first two terms satisfy homogeneous equation (9). Consider the third term
\[D_{0y}^\alpha \left( {\int\limits_0^y {\int\limits_{ - \infty }^{ + \infty } {f\left( {\xi ,\eta } \right){\Gamma _{\alpha  - \frac{\alpha }{{2n}} - 1}}\left( {x - \xi ,y - \eta } \right)} } d\xi d\eta } \right) = \]
\[ = \frac{1}{{\Gamma \left( {2 - \alpha } \right)}}\frac{{{\partial ^2}}}{{\partial {y^2}}}\left( {\int\limits_{ - \infty }^{ + \infty } {d\xi } D_{0y}^{\alpha  - 2}\left( {\int\limits_0^y {f\left( {\xi ,\eta } \right){\Gamma _{\alpha  - \frac{\alpha }{{2n}} - 1}}\left( {x - \xi ,y - \eta } \right)} d\eta } \right)} \right) \Rightarrow \]
\[D_{0y}^{\alpha  - 2}\left( {\int\limits_0^y {f\left( {\xi ,\eta } \right){\Gamma _{\alpha  - \frac{\alpha }{{2n}} - 1}}\left( {x - \xi ,y - \eta } \right)} d\eta } \right) = \]
\[ = \int\limits_0^y {\frac{{d\tau }}{{{{\left( {y - \tau } \right)}^{\alpha  - 1}}}}\left( {\int\limits_0^\tau  {f\left( {\xi ,\eta } \right){\Gamma _{\alpha  - \frac{\alpha }{{2n}} - 1}}\left( {x - \xi ,\tau  - \eta } \right)} d\eta } \right)}  = \]
\[ = \int\limits_0^y {f\left( {\xi ,\eta } \right)d\eta \left( {\int\limits_\eta ^y {\frac{{{\Gamma _{\alpha  - \frac{\alpha }{{2n}} - 1}}\left( {x - \xi ,\tau  - \eta } \right)}}{{{{\left( {y - \tau } \right)}^{\alpha  - 1}}}}} d\tau } \right)}=\]
\[ = \int\limits_0^y {f\left( {\xi ,\eta } \right)d\eta \left( {\int\limits_0^1 {\frac{{{{\left( {y - \eta } \right)}^{2 - \alpha }}{\Gamma _{\alpha  - \frac{\alpha }{{2n}} - 1}}\left( {x - \xi ,\left( {y - \eta } \right)t} \right)}}{{{{\left( {1 - t} \right)}^{\alpha  - 1}}}}} dt} \right)}  = \]
\[ = \Gamma \left( {2 - \alpha } \right)\int\limits_0^y {f\left( {\xi ,\eta } \right){\Gamma _{1 - \frac{\alpha }{{2n}}}}\left( {x - \xi ,y - \eta } \right)d\eta }.\eqno(21)\]
Hence we have
\[D_{0y}^\alpha \int\limits_{ - \infty }^{ + \infty } {d\xi } \left( {\int\limits_0^y {f\left( {\xi ,\eta } \right){\Gamma _{\alpha  - \frac{\alpha }{{2n}} - 1}}\left( {x - \xi ,y - \eta } \right)} d\eta } \right) = \]
\[ = \int\limits_{ - \infty }^{ + \infty } {d\xi } \left( {\int\limits_0^y {f\left( {\xi ,\eta } \right){\Gamma _{ - \frac{\alpha }{{2n}} - 1}}\left( {x - \xi ,y - \eta } \right)} d\eta } \right).\eqno(22)\]
Note that taking into account (21), we have
\[\mathop {\lim }\limits_{y \to  + 0} D_{0y}^{\alpha  - 2}\left( {\int\limits_{ - \infty }^{ + \infty } {d\xi } \left( {\int\limits_0^y {f\left( {\xi ,\eta } \right){\Gamma _{\alpha  - 1 - \frac{\alpha }{{2n}}}}\left( {x - \xi ,y - \eta } \right)} d\eta } \right)} \right) = \]
\[ = \mathop {\lim }\limits_{y \to  + 0} \int\limits_{ - \infty }^{ + \infty } {d\xi } \left( {\int\limits_0^y {f\left( {\xi ,\eta } \right){\Gamma _{1 - \frac{\alpha }{{2n}}}}\left( {x - \xi ,y - \eta } \right)} d\eta } \right) = 0,\eqno(23)\]
\[\mathop {\lim }\limits_{y \to  + 0} D_{0y}^{\alpha  - 1}\left( {\int\limits_{ - \infty }^{ + \infty } {d\xi } \left( {\int\limits_0^y {f\left( {\xi ,\eta } \right){\Gamma _{\alpha  - 1 - \frac{\alpha }{{2n}}}}\left( {x - \xi ,y - \eta } \right)} d\eta } \right)} \right) = \]
\[ = \mathop {\lim }\limits_{y \to  + 0} \int\limits_{ - \infty }^{ + \infty } {d\xi } \left( {\int\limits_0^y {f\left( {\xi ,\eta } \right){\Gamma _{ - \frac{\alpha }{{2n}}}}\left( {x - \xi ,y - \eta } \right)} d\eta } \right) = 0.\eqno(24)\]
From (23), (24) and Lemma 3 we have
\[\begin{array}{l}
\mathop {\lim }\limits_{y \to  + 0} D_{0y}^{\alpha  - 2}u\left( {x,y} \right) = \varphi \left( x \right),\\
\mathop {\lim }\limits_{y \to  + 0} D_{0y}^{\alpha  - 1}u\left( {x,y} \right) = \psi \left( x \right).
\end{array}\]
Now we find the partial derivatives with respect to the variable $ x $, the third term
\[\frac{\partial }{{\partial x}}\left( {\int\limits_{ - \infty }^x {d\xi } \left( {\int\limits_0^y {f\left( {\xi ,\eta } \right){\Gamma _{\alpha  - \frac{\alpha }{{2n}} - 1}}\left( {x - \xi ,y - \eta } \right)} d\eta } \right)} \right) + \]
\[ + \frac{\partial }{{\partial x}}\left( {\int\limits_x^{ + \infty } {d\xi } \left( {\int\limits_0^y {f\left( {\xi ,\eta } \right){\Gamma _{\alpha  - \frac{\alpha }{{2n}} - 1}}\left( {x - \xi ,y - \eta } \right)} d\eta } \right)} \right) = \]
\[ = \mathop {\lim }\limits_{\xi  \to x} \int\limits_0^y {f\left( {\xi ,\eta } \right)\left( {\Gamma _{\alpha  - 1 - \frac{\alpha }{{2n}}}^1\left( {x - \xi ,y - \eta } \right) - \Gamma _{\alpha  - 1 - \frac{\alpha }{{2n}}}^2\left( {x - \xi ,y - \eta } \right)} \right)} d\eta  + \]
\[ + \int\limits_{ - \infty }^{ + \infty } {d\xi } \left( {\int\limits_0^y {f\left( {\xi ,\eta } \right)\frac{\partial }{{\partial x}}{\Gamma _{\alpha  - \frac{\alpha }{{2n}} - 1}}\left( {x - \xi ,y - \eta } \right)} d\eta } \right) = \]
\[ = \int\limits_{ - \infty }^{ + \infty } {d\xi } \left( {\int\limits_0^y {f\left( {\xi ,\eta } \right)\frac{\partial }{{\partial x}}{\Gamma _{\alpha  - \frac{\alpha }{{2n}} - 1}}\left( {x - \xi ,y - \eta } \right)} d\eta } \right),\]
continuing this process we get
\[\frac{{{\partial ^{2n}}}}{{\partial {x^{2n}}}}\left( {\int\limits_{ - \infty }^{ + \infty } {d\xi } \left( {\int\limits_0^y {f\left( {\xi ,\eta } \right){\Gamma _{\alpha  - \frac{\alpha }{{2n}} - 1}}\left( {x - \xi ,y - \eta } \right)} d\eta } \right)} \right) = \]
\[ = \mathop {\lim }\limits_{\varepsilon  \to  + 0} \left( {\frac{\partial }{{\partial x}}\left( {\int\limits_{ - \infty }^{x - \varepsilon } {d\xi } \left( {\int\limits_0^y {f\left( {\xi ,\eta } \right)\frac{{{\partial ^{2n - 1}}}}{{\partial {x^{2n - 1}}}}{\Gamma _{\alpha  - \frac{\alpha }{{2n}} - 1}}\left( {x - \xi ,y - \eta } \right)} d\eta } \right)} \right)} \right) + \]
\[ + \mathop {\lim }\limits_{\varepsilon  \to  + 0} \left( {\frac{\partial }{{\partial x}}\left( {\int\limits_{x - \varepsilon }^{x + \varepsilon } {d\xi } \left( {\int\limits_0^y {f\left( {\xi ,\eta } \right)\frac{{{\partial ^{2n - 1}}}}{{\partial {x^{2n - 1}}}}{\Gamma _{\alpha  - \frac{\alpha }{{2n}} - 1}}\left( {x - \xi ,y - \eta } \right)} d\eta } \right)} \right)} \right) + \]
\[ + \mathop {\lim }\limits_{\varepsilon  \to  + 0} \left( {\frac{\partial }{{\partial x}}\left( {\int\limits_{x + \varepsilon }^{ + \infty } {d\xi } \left( {\int\limits_0^y {f\left( {\xi ,\eta } \right)\frac{{{\partial ^{2n - 1}}}}{{\partial {x^{2n - 1}}}}{\Gamma _{\alpha  - \frac{\alpha }{{2n}} - 1}}\left( {x - \xi ,y - \eta } \right)} d\eta } \right)} \right)} \right) = \]
\[ = {I_1} + {I_2} + {I_3},\]
We calculate each term separately
\[{I_1} = \mathop {\lim }\limits_{\varepsilon  \to  + 0} \left( {\int\limits_0^y {f\left( {x - \varepsilon ,\eta } \right)\frac{{{\partial ^{2n - 1}}}}{{\partial {x^{2n - 1}}}}\Gamma _{\alpha  - 1 - \frac{\alpha }{{2n}}}^1\left( {\varepsilon ,y - \eta } \right)} d\eta } \right) + \]
\[ + \int\limits_{ - \infty }^x {d\xi } \left( {\int\limits_0^y {f\left( {\xi ,\eta } \right)\frac{{{\partial ^{2n}}}}{{\partial {x^{2n}}}}{\Gamma _{\alpha  - \frac{\alpha }{{2n}} - 1}}\left( {x - \xi ,y - \eta } \right)} d\eta } \right),\]
from here
\[\mathop {\lim }\limits_{\varepsilon  \to  + 0} \int\limits_0^y {f\left( {x - \varepsilon ,\eta } \right)\frac{{{\partial ^{2n - 1}}}}{{\partial {x^{2n - 1}}}}\Gamma _{\alpha  - 1 - \frac{\alpha }{{2n}}}^1\left( {\varepsilon ,y - \eta } \right)} d\eta  = \]
\[ = \frac{{{{\left( { - 1} \right)}^{n - 1}}}}{{2n}}\sum\limits_{k = 0}^{n - 1} {\mathop {\lim }\limits_{\varepsilon  \to  + 0} \int\limits_0^y {f\left( {x - \varepsilon ,\eta } \right){{\left( {y - \eta } \right)}^{ - 1}}\sum\limits_{m = 1}^\infty  {\frac{{{{\left( { - {e^{\frac{{n - 1 - 2k}}{{2n}}\pi i}}\varepsilon {{\left( {y - \eta } \right)}^{ - \frac{\alpha }{{2n}}}}} \right)}^k}}}{{k!\Gamma \left( { - \frac{\alpha }{{2n}}k} \right)}}} } d\eta }  = \]
\[\left( {z = \varepsilon {{\left( {y - \eta } \right)}^{ - \frac{\alpha }{{2n}}}} \Rightarrow \eta  = y - {\varepsilon ^{\frac{{2n}}{\alpha }}}{z^{ - \frac{{2n}}{\alpha }}},\,\,d\eta  = \frac{{2n}}{\alpha }{\varepsilon ^{\frac{{2n}}{\alpha }}}{z^{ - \frac{{2n}}{\alpha } - 1}}dz} \right)\]
\[ = \frac{{{{\left( { - 1} \right)}^{n - 1}}}}{\alpha }f\left( {x,y} \right)\sum\limits_{k = 0}^{n - 1} {\sum\limits_{m = 1}^\infty  {\frac{{{{\left( { - {e^{\frac{{n - 1 - 2k}}{{2n}}\pi i}}} \right)}^m}}}{{m!\Gamma \left( { - \frac{\alpha }{{2n}}m} \right)}}\int\limits_0^{ + \infty } {{z^{m - 1}}} } dz}  = \]
\[ = \frac{{{{\left( { - 1} \right)}^n}}}{{2n}}f\left( {x,y} \right)\sum\limits_{k = 0}^{n - 1} {\left. {\left( {\sum\limits_{m = 0}^\infty  {\frac{{{{\left( { - {e^{\frac{{n - 1 - 2k}}{{2n}}\pi i}}} \right)}^m}{z^m}}}{{m!\Gamma \left( { - \frac{\alpha }{{2n}}m + 1} \right)}}}  - 1} \right)} \right|_{z = 0}^{z = \infty }}  = \frac{{{{\left( { - 1} \right)}^{n - 1}}}}{2}f\left( {x,y} \right).\]
Means
\[{I_1} = \frac{{{{\left( { - 1} \right)}^{n - 1}}}}{2}f\left( {x,y} \right) + \int\limits_{ - \infty }^x {d\xi } \left( {\int\limits_0^y {f\left( {\xi ,\eta } \right)\frac{{{\partial ^{2n}}}}{{\partial {x^{2n}}}}{\Gamma _{\alpha  - \frac{\alpha }{{2n}} - 1}}\left( {x - \xi ,y - \eta } \right)} d\eta } \right).\]
Similarly
\[{I_3} = \frac{{{{\left( { - 1} \right)}^{n - 1}}}}{2}f\left( {x,y} \right) + \int\limits_x^{ + \infty } {d\xi } \left( {\int\limits_0^y {f\left( {\xi ,\eta } \right)\frac{{{\partial ^{2n}}}}{{\partial {x^{2n}}}}{\Gamma _{\alpha  - \frac{\alpha }{{2n}} - 1}}\left( {x - \xi ,y - \eta } \right)} d\eta } \right),\]
it is easy to show that
\[{I_2} = 0.\]
so
\[\frac{{{\partial ^{2n}}}}{{\partial {x^{2n}}}}\left( {\int\limits_{ - \infty }^{ + \infty } {d\xi } \left( {\int\limits_0^y {f\left( {\xi ,\eta } \right){\Gamma _{\alpha  - \frac{\alpha }{{2n}} - 1}}\left( {x - \xi ,y - \eta } \right)} d\eta } \right)} \right) = \]
\[ = {\left( { - 1} \right)^{n - 1}}f\left( {x,y} \right) + \int\limits_{ - \infty }^{ + \infty } {d\xi } \left( {\int\limits_0^y {f\left( {\xi ,\eta } \right)\frac{{{\partial ^{2n}}}}{{\partial {x^{2n}}}}{\Gamma _{\alpha  - \frac{\alpha }{{2n}} - 1}}\left( {x - \xi ,y - \eta } \right)} d\eta } \right) = \]
\[ = {\left( { - 1} \right)^{n - 1}}f\left( {x,y} \right) + {\left( { - 1} \right)^{n - 1}}\int\limits_{ - \infty }^{ + \infty } {d\xi } \left( {\int\limits_0^y {f\left( {\xi ,\eta } \right){\Gamma _{ - \frac{\alpha }{{2n}} - 1}}\left( {x - \xi ,y - \eta } \right)} d\eta } \right).\eqno(25)\]
From (22) and (25) we obtain that the expression
\[ - \int\limits_0^y {\int\limits_{ - \infty }^{ + \infty } {f\left( {\xi ,\eta } \right){\Gamma _{\alpha  - \frac{\alpha }{{2n}} - 1}}\left( {x - \xi ,y - \eta } \right)} } d\xi d\eta,\]
satisfies equation (9). The theorem is proved. \\
The question of the uniqueness of the Cauchy problem is still open.
Function
\[-{{\Gamma _{\alpha  - \frac{\alpha }{{2n}} - 1}}\left( {x - \xi ,y - \eta } \right)}\eqno(26)\]
will be called the fundamental solution of equation (9). Representation (26) coincides with the form of the fundamental solution obtained for the case $ 0 <\alpha < 1 $ in [5].
\begin{center}
\textbf{3. Application
}\end{center}
In [7], for the equation of beam vibration
\[{u_{tt}}\left( {t,x} \right) + {a ^2}{u_{xxxx}}\left( {t,x} \right) = 0,\]
a representation of the solution of the Cauchy problem was obtained using the previously constructed fundamental solution, which for $ a = 1 $ has the form
\[\begin{array}{l}
{G_2}\left( {x,t,\xi } \right) = \sqrt {\frac{t}{\pi }} \sin \left[ {\frac{{{{\left( {\xi  - x} \right)}^2}}}{{4t}} + \frac{\pi }{4}} \right] + \\
 + \frac{{\xi  - x}}{2}\left\{ {S\left[ {\left( {\frac{{{{\left( {\xi  - x} \right)}^2}}}{{4t }}} \right)} \right] - C\left[ {\left( {\frac{{{{\left( {\xi  - x} \right)}^2}}}{{4t }}} \right)} \right]} \right\},
\end{array}\]
where
\[S\left( z \right) = \frac{1}{{2\sqrt \pi  }}\int\limits_0^z {\frac{{\sin t}}{{\sqrt t }}} dt,\,\,C\left( z \right) = \frac{1}{{2\sqrt \pi  }}\int\limits_0^z {\frac{{cost}}{{\sqrt t }}} dt - \]
Fresnel integrals.\\
Our fundamental solution (27), for $\alpha=2, \, n=2 $, has the form ($ y $ was replaced by $ t $, $\eta = 0 $)
\[ - {\Gamma _{\frac{1}{2}}}\left( {x - \xi ,t} \right) =  - \frac{{\sqrt t }}{4}\left( { - {e^{\frac{\pi }{4}i}}\sum\limits_{n = 0}^\infty  {\frac{{{{\left( { - {e^{\frac{\pi }{4}i}}} \right)}^n}{\tau ^n}}}{{n!\Gamma \left( { - \frac{1}{2}n + \frac{3}{2}} \right)}}}  - {e^{ - \frac{\pi }{4}i}}\sum\limits_{n = 0}^\infty  {\frac{{{{\left( { - {e^{ - \frac{\pi }{4}i}}} \right)}^n}{\tau ^n}}}{{n!\Gamma \left( { - \frac{1}{2}n + \frac{3}{2}} \right)}}} } \right),\]
where
\[\tau  = \frac{{\left| {x - \xi } \right|}}{{\sqrt t }},\]
show that
$${G_2}\left( {x,t,\xi } \right) =- {\Gamma _{\frac{1}{2}}}\left( {x - \xi ,t} \right).$$
Indeed, we have
\[\sum\limits_{n = 0}^\infty  {\frac{{{{\left( { - {e^{\frac{\pi }{4}i}}} \right)}^n}{\tau ^n}}}{{n!\Gamma \left( { - \frac{1}{2}n + \frac{3}{2}} \right)}}}  = \sum\limits_{n = 0}^\infty  {\frac{{{{\left( { - {e^{\frac{\pi }{4}i}}} \right)}^{4n}}{\tau ^{4n}}}}{{\left( {4n} \right)!\Gamma \left( { - 2n + \frac{3}{2}} \right)}}}  + \sum\limits_{n = 0}^\infty  {\frac{{{{\left( { - {e^{\frac{\pi }{4}i}}} \right)}^{4n + 1}}{\tau ^{4n + 1}}}}{{\left( {4n + 1} \right)!\Gamma \left( { - 2n + 1} \right)}}}  + \]
\[ + \sum\limits_{n = 0}^\infty  {\frac{{{{\left( { - {e^{\frac{\pi }{4}i}}} \right)}^{4n + 2}}{\tau ^{4n + 2}}}}{{\left( {4n + 2} \right)!\Gamma \left( { - 2n + \frac{1}{2}} \right)}}}  + \sum\limits_{n = 0}^\infty  {\frac{{{{\left( { - {e^{\frac{\pi }{4}i}}} \right)}^{4n + 3}}{\tau ^{4n + 3}}}}{{\left( {4n + 3} \right)!\Gamma \left( { - 2n} \right)}}}  = \]
\[ = \sum\limits_{n = 0}^\infty  {\frac{{{{\left( { - 1} \right)}^n}{\tau ^{4n}}}}{{\left( {4n} \right)!\Gamma \left( { - 2n + \frac{3}{2}} \right)}}}  + i\sum\limits_{n = 0}^\infty  {\frac{{{{\left( { - 1} \right)}^n}{\tau ^{4n + 2}}}}{{\left( {4n + 2} \right)!\Gamma \left( { - 2n + \frac{1}{2}} \right)}}}  + \left( { - {e^{\frac{\pi }{4}i}}} \right)\tau  = {I_1}.\]
Similarly
\[\sum\limits_{n = 0}^\infty  {\frac{{{{\left( { - {e^{ - \frac{\pi }{4}i}}} \right)}^n}{\tau ^n}}}{{n!\Gamma \left( { - \frac{1}{2}n + \frac{3}{2}} \right)}}}  = \sum\limits_{n = 0}^\infty  {\frac{{{{\left( { - 1} \right)}^n}{\tau ^{4n}}}}{{\left( {4n} \right)!\Gamma \left( { - 2n + \frac{3}{2}} \right)}}}  - i\sum\limits_{n = 0}^\infty  {\frac{{{{\left( { - 1} \right)}^n}{\tau ^{4n + 2}}}}{{\left( {4n + 2} \right)!\Gamma \left( { - 2n + \frac{1}{2}} \right)}}}  - \tau {e^{ - \frac{\pi }{4}i}} = {I_2},\]
from here
\[ - {\Gamma _{\frac{1}{2}}}\left( {x - \xi ,t} \right) =  - \frac{{\sqrt t }}{4}\left( { - {e^{\frac{\pi }{4}i}}{I_1} - {e^{ - \frac{\pi }{4}i}}{I_2}} \right) = \]
\[ = \frac{{\sqrt {2t} }}{8}\left( {\left( {{I_1} + {I_2}} \right) + i\left( {{I_1} - {I_2}} \right)} \right),\]
further
\[{I_1} + {I_2} = 2\sum\limits_{n = 0}^\infty  {\frac{{{{\left( { - 1} \right)}^n}{\tau ^{4n}}}}{{\left( {4n} \right)!\Gamma \left( { - 2n + \frac{3}{2}} \right)}}}  - \sqrt 2 \tau ,\]
\[\left( {{I_1} - {I_2}} \right)i =  - 2\sum\limits_{n = 0}^\infty  {\frac{{{{\left( { - 1} \right)}^n}{\tau ^{4n + 2}}}}{{\left( {4n + 2} \right)!\Gamma \left( { - 2n + \frac{1}{2}} \right)}}}  + \sqrt 2 \tau  \Rightarrow \]
\[ - {\Gamma _{\frac{1}{2}}}\left( {x - \xi ,t} \right) = \frac{{\sqrt {2t} }}{4}\left( {\sum\limits_{n = 0}^\infty  {\frac{{{{\left( { - 1} \right)}^n}{\tau ^{4n}}}}{{\left( {4n} \right)!\Gamma \left( { - 2n + \frac{3}{2}} \right)}}}  - \sum\limits_{n = 0}^\infty  {\frac{{{{\left( { - 1} \right)}^n}{\tau ^{4n + 2}}}}{{\left( {4n + 2} \right)!\Gamma \left( { - 2n + \frac{1}{2}} \right)}}} } \right),\]
further, calculations show that
\[\sum\limits_{n = 0}^\infty  {\frac{{{{\left( { - 1} \right)}^n}{\tau ^{4n}}}}{{\left( {4n} \right)!\Gamma \left( { - 2n + \frac{3}{2}} \right)}}}  = \frac{{\Gamma \left( { - \frac{1}{2}} \right)}}{\pi }{}_1{F_2}\left[ {\begin{array}{*{20}{c}}
{ - \frac{1}{4}, - {{\left( {\frac{{{\tau ^2}}}{8}} \right)}^2}}\\
{\frac{3}{4},\frac{1}{2}}
\end{array}} \right],\]
\[\sum\limits_{n = 0}^\infty  {\frac{{{{\left( { - 1} \right)}^n}{\tau ^{4n + 2}}}}{{\left( {4n + 2} \right)!\Gamma \left( { - 2n + \frac{1}{2}} \right)}}}  = \frac{{\Gamma \left( { - \frac{1}{2}} \right){\tau ^2}}}{{4\pi }}{}_1{F_2}\left[ {\begin{array}{*{20}{c}}
{\frac{1}{4}, - {{\left( {\frac{{{\tau ^2}}}{8}} \right)}^2}}\\
{\frac{3}{2},\frac{5}{4}}
\end{array}} \right],\]
where
\[{}_1{F_2}\left[ {\begin{array}{*{20}{c}}
{a,z}\\
{b,c}
\end{array}} \right] = \sum\limits_{n = 0}^\infty  {\frac{{{{\left( a \right)}_n}}}{{{{\left( b \right)}_n}{{\left( c \right)}_n}}}\frac{{{z^n}}}{{n!}}} \]
-hypergeometric function,
\[{\left( d \right)_n} = d\left( {d + 1} \right)...\left( {d + n - 1} \right),\]
from here
\[ - {\Gamma _{\frac{1}{2}}}\left( {x - \xi ,t} \right) = \frac{{\sqrt {2t} \Gamma \left( { - \frac{1}{2}} \right)}}{{16\pi }}\left( {4{}_1{F_2}\left[ {\begin{array}{*{20}{c}}
{ - \frac{1}{4}, - {{\left( {\frac{{{\tau ^2}}}{8}} \right)}^2}}\\
{\frac{3}{4},\frac{1}{2}}
\end{array}} \right] - {\tau ^2}{}_1{F_2}\left[ {\begin{array}{*{20}{c}}
{\frac{1}{4}, - {{\left( {\frac{{{\tau ^2}}}{{8t}}} \right)}^2}}\\
{\frac{3}{2},\frac{5}{4}}
\end{array}} \right]} \right).\]
Now, taking into account the relations (see [8])
\[4{}_1{F_2}\left[ {\begin{array}{*{20}{c}}
{ - \frac{1}{4}, - {{\left( {\frac{{{\tau ^2}}}{8}} \right)}^2}}\\
{\frac{3}{4},\frac{1}{2}}
\end{array}} \right] = 4\cos \left( {\frac{{{\tau ^2}}}{4}} \right) + 2\sqrt {2\pi } \left| \tau  \right|S\left( {\frac{{{\tau ^2}}}{4}} \right),\]
\[{\tau ^2}{}_1{F_2}\left[ {\begin{array}{*{20}{c}}
{\frac{1}{4}, - {{\left( {\frac{{{\tau ^2}}}{8}} \right)}^2}}\\
{\frac{3}{2},\frac{5}{4}}
\end{array}} \right] = 2\left| \tau  \right|\sqrt {2\pi } C\left( {\frac{{{\tau ^2}}}{4}} \right) - 4\sin \frac{{{\tau ^2}}}{4},\]
finally get
\[ - {\Gamma _{\frac{1}{2}}}\left( {x - \xi ,t} \right) = \sqrt {\frac{t}{\pi }} \sin \left( {\frac{{{\tau ^2}}}{4} + \frac{\pi }{4}} \right) + \frac{{\left| {\xi  - x} \right|}}{2}\left( {S\left( {\frac{{{\tau ^2}}}{4}} \right) - C\left( {\frac{{{\tau ^2}}}{4}} \right)} \right),\]
i.e.
$$- {\Gamma _{\frac{1}{2}}}\left( {x - \xi ,t} \right)={G_2}\left( {x,t,\xi } \right).$$

\begin{center}
\textbf{References}
\end{center}
 1. Irgashev B. Y. Construction of Singular Particular Solutions Expressed via Hypergeometric Functions for an Equation with Multiple Characteristics. //Differential Equations. 2020. Т. 56. No 3, pp. 315-323.\\
 2.  L u c h k o Yu.,  G o r e n f l o R. Scale-invariant solutions of a partial
differential equation of fractional order. // Fractional Calculus and Applied
Analysis.1998. Т.1, No 1. pp. 63-78. \\
3.  Wright E. M. The generalized Bessel function of order greater than one // Quart.
J. Math., Oxford Ser., 1940. Т. 1, No 11. pp. 36–48.\\
4. Pskhu A. Fundamental solutions and cauchy problems for an odd-order partial differential equation with. // Electronic Journal of Differential Equations, Vol. 2019 (2019), No. 21, pp. 1–13. \\
5. Karasheva L. L. Cauchy problem for high even order parabolic equation with time fractional derivative //Sibirskie Elektronnye Matematicheskie Izvestiya [Siberian Electronic Mathematical Reports].2018. Т. 15. pp. 696-706. [In Russian].\\
6. Karasheva L. L. On Properties of an Integer Function That Generalizes the Wright Function //Itogi Nauki i Tekhniki. Seriya" Sovremennaya Matematika i ee Prilozheniya. Tematicheskie Obzory". 2018. Т. 149. pp. 38-43.[In Russian]\\
7. Sabitov K. B. Cauchy problem for the beam vibration equation //Differential Equations. 2017. Т. 53. No. 5. pp. 658-664.\\
8. Dunaev A.S., Shlychkov V.I. Hypergeometric functions. Ural Federal University. Yekaterinburg. 2017. P.1273 [in Russian].

\end{document}